\documentclass[12pt]{amsart}

\usepackage{amssymb}
\usepackage{mathrsfs}
\usepackage{epic}
\usepackage{eepic}

\newtheorem{thm}{Theorem}[section]
\newtheorem{lem}[thm]{Lemma}
\newtheorem{claim}[thm]{Claim}
\newtheorem{question}[thm]{Question}

\theoremstyle{remark}
\newtheorem{remark}[thm]{Remark}

\theoremstyle{definition}
\newtheorem{defn}[thm]{Definition}

\newcommand\Fol{\textrm{F\o l}}
\newcommand{\symdif}{\triangle}
\newcommand\Nbb{\mathbb{N}}
\newcommand\Zbb{\mathbb{Z}}
\newcommand\Tcal{\mathscr{T}}
\newcommand\Xcal{\mathscr{X}}
\newcommand\Ycal{\mathscr{Y}}
\newcommand\Acal{\mathscr{A}}
\newcommand\tie{\char'136{}}
\newcommand\lex{\mathrm{lex}}
\newcommand{\Seq}[1]{\langle #1 \rangle}
\newcommand{\Th}{{}^{\textrm{th}}}
\newcommand{\seq}[1]{\mathtt{#1}}
\newcommand\Ecal{\mathscr{E}}
\newcommand\Rcal{\mathscr{R}}

\title[Fast growth]
{Fast growth in the F\o lner function for Thompson's group $F$}

\keywords{F\o lner function, tower function, Thompson's group, amenable}

\subjclass[2000]{20F65, 43A07}

\thanks{
I would like to thank Matt Brin
for his careful reading of drafts
of this paper, catching errors,
and suggesting improvements.
In particular, the current formulation
of Lemma \ref{marginal} was suggested by him
(the original formulation provided a weaker estimate).
I would also like to thank the anonymous referee for their
very careful reading and helpful comments and suggestions.
This research was supported in part by
NSF grant DMS--0757507.
}

\author{Justin Tatch Moore}

\address{Department of Mathematics \\ Cornell University\\
Ithaca, NY 14853--4201 \\ USA}

\email{{\tt justin@math.cornell.edu}}
 
\begin{document}

\begin{abstract}
The purpose of this note is to prove a lower bound on the
growth of F\o lner functions for Thompson's group $F$.
Specifically I will prove that, for any
finite generating set $\Gamma \subseteq F$,
there is a constant $C$ such that
$\Fol_{F,\Gamma} (C^n) \geq \exp_{n}(0)$.
\end{abstract}

\maketitle

\section{Introduction}

In this paper we will study the F\o lner function for Thompson's group $F$.
Recall that a finite subset
$A$ of a finitely generated group $G$ is \emph{$\varepsilon$-F\o lner} (with
respect to a finite generating set $\Gamma \subseteq G$) if 
\[
\sum_{\gamma \in \Gamma} |(A \cdot \gamma) \ \symdif A| < \varepsilon |A|
\]
where $\symdif$ denotes symmetric difference.
The \emph{F\o lner function} of $G$ (with respect to $\Gamma$) is defined by
\[
\Fol_{G,\Gamma}(n) = \min \{|A| : A \subseteq G \textrm{ is } \frac{1}{n}\textrm{-F\o lner w.r.t. } \Gamma\}
\]
with $\Fol_{G,\Gamma} (n) = \infty$ if there is no $1/n$-F\o lner set with respect to $\Gamma$.
By F\o lner's criterion, a group
$G$ is amenable if and only if its F\o lner function
(with respect to any finite generating set $\Gamma$) is finite valued.

Thompson's group $F$ has many equivalent definitions;
we will use the formulation in terms of tree diagrams defined below.
The standard presentation of $F$ is infinite, with generators $x_i$ $(i \in \Nbb)$
satisfying $x_i^{-1} x_n x_i = x_{n+1}$ for all $i < n$.
It is well known, however, that $F$ admits the finite presentation
\[
\Seq{A,B \ |\ [AB^{-1},A^{-1}BA] = [AB^{-1},A^{-2}BA^2]= \mathrm{id}}
\]
(see \cite{CFP}).
Geoghegan conjectured that $F$ is not amenable \cite[p. 549]{comb_grp_top}
and at present this problem remains open.\footnote{
In fact Richard Thompson himself had studied the question of the amenablity
of $F$ already by the early 1970s \cite{thompson_letter},
although the question
was not well known until it was independently considered and popularized
by Geoghegan.}
The goal of this paper is to establish the following lower bound on the
F\o lner function for $F$.
\begin{thm} \label{tower_growth}
For every finite symmetric generating set $\Gamma \subseteq F$
there is a constant $C>1$ such that if $A \subseteq F$
is a $C^{-n}$-F\o lner set with 
respect to $\Gamma$, then $A$ contains at least $\exp_{n} (0)$ elements.
In particular $\Fol_{F,\Gamma}$ is not eventually dominated $\exp_p(n)$ for any
finite $p$.
\end{thm}
Here $\exp_p(n)$ is the $p$-fold composition of the exponential function
defined by $\exp_{0} (n) = n$ and $\exp_{p+1}(n) = 2^{\exp_p (n)}$.
If it turns out that $F$ is amenable, then Theorem \ref{tower_growth}
would be a step toward answering (negatively)
the following question of Gromov \cite[p. 578]{entropy_isoperimetry}.

\begin{question}
Is there a primitive recursive function which eventually dominates every 
F\o lner function of an amenable finitely presented group?
\end{question}
While it is known that the F\o lner functions of amenable finitely generated groups can grow arbitrarily
fast \cite{pw_auto_groups}, this is not the case for finitely presented groups
(since there are only countably many such groups).
See \cite{iso_profiles_fg} 
for what can be accomplished via wreath products of $\Zbb$.

This note is organized as follows.
In Section \ref{background},
I will review some of the basic definitions associated with $F$ and
fix some notational conventions.
In Section \ref{marginal_section}, I will introduce the notion of a \emph{marginal set}
and prove some basic lemmas about them.
These are sets which must have small intersections with F\o lner sets.
They play a central role in the proof of the main result of the paper.
Section \ref{act_on_Tcal} recasts the amenability of $F$ in terms of its partial right
action on the finite rooted ordered binary trees $\Tcal$.
Section \ref{partial_section}
defines an operation on elements of $\Tcal$ which exponentially decreases their size
and commutes with the partial right action of $F$.
It is shown that the trees which are trivialized by this operation are
marginal and it is this that allows the proof of Theorem \ref{tower_growth}.

\section{Notation and background}

\label{background}

I will use \cite{CFP} as a general reference for Thompson's group $F$, although the reader is warned
that the notation in the present paper will differ somewhat from that of \cite{CFP}.
Let $\Tcal$ denote the collection of all \emph{finite rooted binary trees}.
For concreteness, we will view elements $T$ of $\Tcal$ as finite sets of
binary sequences which have the following property:
for every infinite binary sequence $x$, there is a unique element of $T$ which is an initial
part of $x$.
Thus an element $T$ of $\Tcal$ is a record of
the addresses of the leaves of the tree which it represents.
The \emph{trivial tree} is the set which consists only of the sequence of length $0$.

The collection of finite binary sequences is equipped with the operation of
concatenation (denoted $u \tie v$), the partial order $\subseteq$ of extension (defined
by $u \subseteq u \tie v$), and the lexicographic order (denoted $u <_\lex v$).
Note that \emph{$u$ extends $v$} includes the possibility that $u = v$.
I will use $\Seq{}$ to denote the sequence of length $0$.
$\Tcal$ can also be characterized as being those finite sets $T$ of binary sequences such that:
\begin{itemize}

\item
no element of $T$ is an initial part of another element of $T$ and

\item 
if $u$ is a binary sequence, then $u \tie \seq{0}$ has an extension in $T$ if and only
if $u \tie \seq{1}$ has an extension in $T$.

\end{itemize}
If $U$ and $V$ are in $\Tcal$, then we will say that $U$ is \emph{dominated by} $V$
if every element of $U$ has a extension in $V$.
If $T$ is in $\Tcal$ and $u$ a finite binary sequence,
define
\[
T/u = \{s : u \tie s \in T\}.
\]
If this set is non-empty, then it is again a member of $\Tcal$ 
(in which case $T/u$ is the tree of descendants of $u$).
Elements of $\Tcal$ come with a canonical ordering provided by $<_\lex$
and phrases such
as \emph{the $i\Th$ element} and \emph{the minimum element}
will always refer to the $<_\lex$-order in this context.
In this paper $\Nbb$ contains $0$ and in particular counting will always start at 0.
The letters $i$, $j$, $k$, and $l$ will always be used to denote natural numbers.

A \emph{tree diagram} is a pair $(L,R)$ of elements of $\Tcal$ such that
$|L| = |R|$.
We view a tree diagram as describing a map of sequences defined by
$$s_i \tie x \mapsto t_i \tie x$$
where $s_i$ and $t_i$ are the $i\Th$ elements of $L$ and $R$ respectively and
$x$ is any binary sequence.
This map is defined not only on all infinite length binary sequences but also on all but finitely many
finite binary sequences.
The value of the associated map $f$ at a sequence $t$ will be denoted $t \cdot f$ when it
is defined.
If $t$ is a finite sequence, then we will say that $f$
\emph{acts properly} on $t$ if $t \cdot f$ is defined and
the final digit of $t \cdot f$ agrees with that of $t$.
Notice that $f$ acts properly on $t$ unless $s\cdot f$ is undefined for every
proper initial part $s$ of $t$.
In particular, if $f$ acts properly on $t$, it acts properly on any extension of $t$.

If two tree diagrams define the same map on infinite sequences, then they are said
to be equivalent.
Every tree diagram is equivalent to
a unique minimal tree diagram;
such a tree diagram is said to be \emph{reduced}.
Furthermore a tree diagram $(S,T)$ is reduced if and only if whenever
$i < |S| - 1$ it is not the case that both $s_i$ and $t_i$ end with $\seq{0}$ and
$s_{i+1}$ and $t_{i+1}$ end with $\seq{1}$ (where $s_i$ and $t_i$ are the $i\Th$ elements
of $S$ and $T$ respectively).
See \cite{CFP} for details.

Thompson's group $F$ is the collection of reduced tree diagrams with
the operation defined by $f \cdot g = g \circ f$ (i.e. ``$f$ followed by $g$'').
$F$ is generated by $\{x_0,x_1\}$ where $x_0$ and $x_1$ are specified by:
\[
x_0 =
\begin{cases}
\seq{00} \mapsto \seq{0} \\
\seq{01} \mapsto \seq{10} \\
\seq{1} \mapsto \seq{11} \\
\end{cases}
x_1 =
\begin{cases}
\seq{0} \mapsto \seq{0} \\
\seq{100} \mapsto \seq{10} \\
\seq{101} \mapsto \seq{110} \\
\seq{11} \mapsto \seq{111} \\
\end{cases}
\]
In our discussion of $F$
``generator'' will mean an element of the set $\Gamma = \{x_0,x_1,x_0^{-1},x_1^{-1}\}$
(this is really only relevant in Lemma \ref{doubling}).
This includes the usage of \emph{F\o lner} in Sections \ref{act_on_Tcal} and \ref{partial_section}.
If $x_n$ is defined by $x_{n+1} = x_0^{-n} x_1 x_0^n$ for $n \geq 1$, then we obtain
the generators which yield the infinite presentation of $F$ mentioned in the introduction.

Notice that the existence of constant $C$ satisfying the conclusion of Theorem \ref{tower_growth}
for $\Gamma$ implies the main theorem for all finite generating sets.
This is because if $\Gamma'$ is any other finite generating set, there is a constant $K >0$ such
that any set which is $\varepsilon$-F\o lner with respect to $\Gamma'$ is $K\varepsilon$-F\o lner
with respect to $\Gamma$.

I will also identify elements of $F$ with the corresponding maps on sequences.
If $T$ is in $\Tcal$ and $f$ is in $F$ with $f$ defined on all of $T$,
then $T \cdot f$ is the pointwise image of $T$ under $f$.
It is easily checked that this results in an element of $\Tcal$ and hence
this defines a partial right action of $F$ on $\Tcal$.
If $T$ is in $\Tcal$ and $f$ is in $F$, then $f$ \emph{acts properly} on $T$ if
it acts properly on the elements of $T$.
If $\Gamma$ is a subset of $F$, then $\Gamma$ \emph{acts properly on} $T$ if each element of
$\Gamma$ acts properly on $T$.
Observe that if $(S,T)$ is a reduced tree diagram which represents $g$ and $f$ acts properly
on $T$, then $(S,T \cdot f)$ is reduced and represents $f \circ g$.
If $f$ is in $F$, I will write $(L_f,R_f)$ to denote the reduced tree diagram
for $f$.

\section{Marginal sets}

\label{marginal_section}
In this section we will introduce the notion of a \emph{marginal set} and
collect some basic lemmas which we will use in Sections \ref{act_on_Tcal} and \ref{partial_section}.
Throughout this section ``right'' in an implicit adjective whenever
applicable, although all statements have their corresponding ``left'' analogs.
Fix, for the duration of the section,
a group $G$ with a finite symmetric generating set $\Gamma$.

\begin{defn}
A \emph{partial action} of $G$ on a set $S$ is a partial function
$\cdot: S \times G \to S$ such that:
\begin{itemize}

\item $x \cdot e = x$ for all $x \in S$;

\item $x \cdot g = y$ if and only if $x = y \cdot g^{-1}$ for all $g \in G$ and $x,y \in S$;

\item
$x \cdot (gh) = (x \cdot g) \cdot h$ for all $g,h\in G$ and all $x \in S$ for which
all computations involving $\cdot$ are defined.

\end{itemize}
If $E \subseteq S$ and $g$ is in $G$, I will write $E \cdot g$ to denote $\{x \cdot g : x \in E\}$.
\end{defn}

\begin{remark}
Since the exponentiation is defined in $G$, it is entirely possible in
general for $x \cdot g^i$ to be undefined while $x \cdot g^j$ is defined for
some $i < j$.
\end{remark}

\begin{defn}
A \emph{weighted $\varepsilon$-F\o lner set with respect to $\Gamma$} is a function $\mu$ from a
finite subset of $S$ into $(0,\infty)$
which satisfies
\[
\sum_{\gamma \in \Gamma} \sum_{s \in S} |\mu(s \cdot \gamma) - \mu (s)| < \varepsilon \sum_{s \in S} \mu(s)
\]
where we adopt with the conventions that
$\mu (s) = 0$ if $s$ is not in the domain of $\mu$ and
$\mu(s \cdot g) = 0$ if $s \cdot g$ is undefined.
The set $\{s \in S : \mu(s) > 0\}$ will be referred to as the \emph{support}
of $\mu$.
The function $\mu$ induces a finitely supported measure on $S$, also denoted
$\mu$, defined by $\mu(A) = \sum_{s \in A} \mu(s)$.
In the remainder of the paper, the generating set
will always be clear from the context
and we will often suppress mention of it.
\end{defn}

Notice that an $\varepsilon$-F\o lner set is a set $A \subseteq S$ such that
the characteristic function $1_A$ is a weighted $\varepsilon$-F\o lner set.
I will use $\mu \restriction A$ to denote $\mu \cdot 1_A$.
We will need the following property of weighted F\o lner sets, which is
what justifies their added generality in the present paper.

\begin{lem} \label{Folner_hom_prop}
Suppose $G$ acts partially on sets $S$ and $T$ and
that $\mu$ is a weighted $\varepsilon$-F\o lner set with respect to the action on $S$.
If $h:S \to T$ satisfies that $h(s \cdot \gamma) = h(s) \cdot \gamma$
(with both quantities defined) whenever 
$\mu(s) + \mu(s \cdot \gamma) > 0$ and $\gamma \in \Gamma$,
then $\nu(t) = \sum_{h(s) = t} \mu(s)$ ($s$ ranges over $S$) defines a weighted $\varepsilon$-F\o lner
set (with respect to the action on $T$). 
\end{lem}

\begin{proof}
Let $G$, $S$, $T$, $\mu$, and $\nu$ be as in the statement of the lemma.
We need to verify that
\[
\sum_{\gamma \in \Gamma} \sum_{t \in T} |\nu(t \cdot \gamma) - \nu (t)| < \varepsilon \nu(T)                             
\]
Let $\gamma \in \Gamma$ and $t \in T$ be fixed for the moment.
First observe that if $|\nu(t \cdot \gamma) - \nu (t)| > 0$,
then $t \cdot \gamma$ is defined.
This is because otherwise it must be the case that $\nu(t) > 0$ and
hence there must be an $s$ such that
$h(s) = t$ and $\mu(s) > 0$.
In particular, this implies $t \cdot \gamma$ is defined.

Next we have that
\[
\sum_{h(u) = t \cdot \gamma} \mu(u) = \sum_{h(u) \cdot \gamma^{-1} = t} \mu(u) =
\sum_{h(u \cdot \gamma^{-1}) = t} \mu (u)
= \sum_{h(s) = t} \mu(s \cdot \gamma)
\]
The first equality is justified by the properties of a partial action;
the second equality is justified by our assumption that $\mu(u) > 0$ implies
$h(u \cdot \gamma^{-1}) = h(u) \cdot \gamma^{-1}$ with both quantities defined;
the third equality is justified by the properties of a partial action and our assumption
that $\mu(u) > 0$ implies $u \cdot \gamma^{-1}$ is defined.
Now it follows that
\[
\sum_{\gamma \in \Gamma} \sum_{t \in T} |\nu(t \cdot \gamma) - \nu (t)| =
\sum_{\gamma \in \Gamma} \sum_{t \in T} |\sum_{h(s) = t} \mu(s \cdot \gamma) - \mu(s)|                             
\]
\[
\leq \sum_{\gamma \in \Gamma} \sum_{s \in S} |\mu(s \cdot \gamma) - \mu(s)| < \varepsilon \mu(S) = \varepsilon \nu(T).
\]
\end{proof}

Fix a partial action of $G$ on a set $S$ for the duration of this section.
If $g$ is in $G$, let $d_g$ be the minimum length of a word in $\Gamma$
which evaluates to $g$.
We will need the following lemma.
\begin{lem} \label{dg_Folner}
If $\varepsilon > 0$, $\mu$ is a weighted $\varepsilon$-F\o lner set and 
$g$ is in $G$, then
\[
\sum_{s \in S} |\mu(s \cdot g) - \mu(s)| < 2\varepsilon d_g \mu(S).
\]
\end{lem}

\begin{proof}
Let $\gamma_i$ $(i < d_g)$ be elements of $\Gamma$
such that $g = \prod_{i < d_g} \gamma_i$.
Let $g_j = \prod_{i < j} \gamma_i$.
\[
\sum_{s \in S} |\mu(s \cdot g) - \mu(s)| \leq \sum_{s \in S} \sum_{i < d_g} |\mu(s \cdot g_{i+1}) - \mu(s \cdot g_i)|.
\]
Notice that it may be that $s \cdot g_{i+1}$ is defined even though $s \cdot g_i$ is not.
If this is the case, however, then
\[
|\mu(s \cdot g_{i+1}) - \mu(s \cdot g_i)| = |\mu((s \cdot g_{i+1}) \cdot \gamma_i^{-1}) - \mu(s \cdot g_{i+1})|.
\]
It follows that
\[
\sum_{i < d_g} \sum_{s \in S} |\mu(s \cdot g_{i+1}) - \mu(s \cdot g_i)| \leq
\sum_{i < d_g} \sum_{\sigma = \pm 1} \sum_{s \in S} |\mu(s \cdot \gamma_i^\sigma) - \mu(s)|
\]
which is less than $2\varepsilon d_g \mu (S)$.
\end{proof}

We will be interested in the following notion which ensures that a set
has small intersection with any F\o lner set.
The definition is motivated by the following simple observation.
If $\mu$ is a finitely additive invariant probability measure on a group $G$
and $E \subseteq G$ satisfies that, for some $g \in G$,
$\{E \cdot g^i : i \in \Nbb\}$ is a pairwise disjoint family, then
$\mu(E) = 0$.

\begin{defn}
If $g \in G$, $I \subseteq S$, and $E \subseteq S$, then
\emph{$g$ marginalizes $E$ off $I$} if for every $x \in E$
if $x \cdot g^k \in E$ and $k > 0$, then there is an $i < k$ such that
$x \cdot g^i$ is in $I$ or is undefined.
If $I$ is the emptyset, then I will write $g$ \emph{marginalizes} $E$.
\end{defn}

\begin{defn}
The $k$-\emph{marginal sets} for the partial action of $G$ on $S$ are
defined recursively as follows.
The emptyset is $0$-marginal.
If there is a decomposition
$E = \bigcup_{i < l} E_i \subseteq S$ and for each $i < l$, there is
a $g_i \in G$ and a $k$-marginal set $I_i$ such that
$g_i$ marginalizes $E_i$ off $I_i$, then $E$ is $(k+1)$-marginal.
$E \subseteq S$ is \emph{marginal}
if it is $k$-marginal for some $k < \infty$.
\end{defn}

\begin{remark} \label{marginal_rem}
For a fixed partial action of a group $G$ on a set $S$,
it is immediate from the definition that a finite union of marginal sets
is marginal and that a subset of a marginal set is marginal.
Additionally, if $E \subseteq S$ and $g$ is in $G$, then
$g^{-1}$ marginalizes $E \setminus (E \cdot g)$ off $E$.
In particular, if $E$ is marginal, then so is $E \cdot g$.
\end{remark}

We will need the following lemma which shows that marginal sets have small
intersections with F\o lner sets.

\begin{lem} \label{marginal}
Suppose $\varepsilon > 0$ and
$\mu$ is a weighted $\varepsilon$-F\o lner set with support $A$.
If $E \subseteq S$ and $g$ marginalizes $E$ off $S \setminus A$,
then $\mu(E) < d_g \varepsilon \mu(S)$.
\end{lem}

\begin{proof}
For each $x \in E \cap A$, let
$\Phi(x)$ denote the set of all $x \cdot g^i$ such that
\begin{itemize}

\item
$x \cdot g^j$ is defined and in the support of $\mu$
for all $j \leq i$ and

\item
$\mu(x \cdot g^{i+1}) < \mu (x \cdot g^j)$ for all $j \leq i$.

\end{itemize}
Observe that each $\Phi(x)$ is finite and non empty.
Because $E$ is marginalized off $S \setminus A$ by $g$,
$\{\Phi(x) : x \in E\}$ is a pairwise disjoint family.
Observe that for a fixed $x \in E$
\[
\mu(x) \leq
\sum_{y \in \Phi(x)} \mu(y) - \mu(y \cdot g) =
\sum_{y \in \Phi(x)} |\mu(y \cdot g) - \mu (y)|.
\]
Summing over $x \in E$ and applying Lemma \ref{dg_Folner}
now gives the desired estimate.
\end{proof}

\begin{lem}
If $E \subseteq S$ is marginal, then there is a constant $C$ such that
if $\varepsilon > 0$ and $\mu$ is a weighted $\varepsilon$-F\o lner set, then
$\mu(E) < C \varepsilon \mu(S)$.
\end{lem}

\begin{proof}
I will prove by induction on $k$ that if $E$ is $k$-marginal, then the
conclusion of the lemma holds.
If $k = 0$, then $E$ is in fact empty and $\mu(E) = 0$.
Next suppose that $E$ is $(k+1)$-marginal.
Let $I_i$ $(i < l)$, $g_i$ $(i < l)$, and $E_i$ $(i < l)$ be such that:
\begin{itemize}

\item for each $i < l$, $g_i$ is in $G$ and marginalizes $E_i$ off $I_i$;

\item for each $i < l$, $I_i$ is $k$-marginal;

\item $E = \bigcup_{i < l} E_i$.

\end{itemize}
By our inductive assumption, there are $C_i$ $(i < l)$ such that 
if $\mu$ is a weighted $\varepsilon$-F\o lner set, then
$\mu(I_i) < C_i \varepsilon \mu(S)$.
Set $C = \sum_{i < l} C_i + 2 d_{g_i}$.
Now if $\mu$ is a weighted $\varepsilon$-F\o lner set,
\[
\mu(E) = \mu(\bigcup_{i < l} I_i) + \mu(\bigcup_{i < l} E_i \setminus \bigcup_{i < l} I_i)
\leq \sum_{i < l} C_i \varepsilon \mu(S) + 2 d_{g_i} \varepsilon \mu(S) = C \varepsilon \mu(S).
\]
\end{proof}

\begin{lem} \label{Folner_refine}
If $\mu$ is a weighted $\varepsilon$-F\o lner set and $\nu$ is a function from a finite subset
of $S$ into $[0,\infty)$ with $\nu \leq \mu$ pointwise and $\nu(S) \geq (1-\delta) \mu(S)$ for some
$\delta$ satisfying $0 < \delta < 1$,
then $\nu$ is a weighted $(\varepsilon +2 |\Gamma| \delta)/(1-\delta)$-F\o lner set.
\end{lem}

\begin{proof}
Let $\mu$ and $\nu$ be given as in the statement of the lemma.
Fix an $s \in S$ and $\gamma \in \Gamma$.
If $\nu(s \cdot \gamma) \geq \nu (s)$, then
\[
|\nu(s \cdot \gamma) - \nu(s)| \leq \mu(s \cdot \gamma) - \mu(s) + \mu(s)-\nu(s).
\]
If $\nu(s \cdot \gamma) \leq \nu(s)$, then
\[
|\nu(s \cdot \gamma) - \nu(s)| \leq \mu(s) - \mu(s \cdot \gamma) + \mu(s \cdot \gamma) - \nu(s \cdot \gamma).
\]
In either case, we have
\[
|\nu(s \cdot \gamma) - \nu(s)| \leq |\mu(s \cdot \gamma) - \mu(s)| + \mu(s) - \nu(s) + \mu(s \cdot \gamma) - \nu (s \cdot \gamma).
\]
Summing over $s \in S$ and $\gamma \in \Gamma$ and combining this with our
hypotheses we obtain:
\[
\sum_{\gamma \in \Gamma} \sum_{s \in S} |\nu(s \cdot \gamma) - \nu(s)| < \varepsilon \mu(S) + 2 |\Gamma| \delta \mu(S)
\leq \frac{ \varepsilon + 2|\Gamma| \delta}{1 - \delta} \nu(S).
\]
\end{proof}

\begin{lem} \label{marginal->Folner}
If $E \subseteq S$ is marginal and $S \setminus E$ is non-empty,
then there is a constant $C$ such that if $\mu$ is a weighted $\varepsilon$-F\o lner set
and $C \varepsilon \leq 1$,
then $\mu \restriction (S \setminus E)$ is a weighted $C \varepsilon$-F\o lner (and in particular the support of
$\mu \restriction (S \setminus E)$ is non-empty).
\end{lem}

\begin{proof}
Let $E$ be marginal and let $C_0$ be such that if $\varepsilon > 0$ and $\mu$ is a weighted
$\varepsilon$-F\o lner set, then
$\mu(E) < C_0 \varepsilon \mu(S)$ and hence
$\mu(S \setminus E)  > C_0 \varepsilon \mu (S)$.
Set $C = 2(1 + 2 C_0 |\Gamma|)$ and suppose that $\varepsilon > 0$
satisfies $C \varepsilon \leq 1$ and that $\mu$ is
a weighted $\varepsilon$-F\o lner set.
It follows that $C_0 \varepsilon \leq 1/2$ and hence
\[
\frac{\varepsilon + 2 |\Gamma| C_0 \varepsilon}{1-C_0 \varepsilon} \leq
2(1 + 2 |\Gamma| C_0) \varepsilon.
\]
Observe that the support of $\mu$ is not contained in $E$ since
$\mu(E) < C_0 \varepsilon \mu(S) < \mu (S)$.
Hence by applying Lemma \ref{Folner_refine} with $\delta = C_0 \varepsilon$
and $\nu = \mu \restriction (S \setminus E)$, we obtain that 
$\mu \restriction (S \setminus E)$ is $C \varepsilon$-F\o lner.
\end{proof}

\begin{defn}
A subset $A \subseteq G$ is \emph{$\Gamma$-connected}
if whenever $x$ and $y$ are in $A$,
there are $\gamma_i$ $(i < l)$ in $\Gamma$ such that, setting
$x_0 = x$ and $x_{i+1} = x_i \cdot \gamma_i$, then $x_i$ is defined
for each $i \leq l$ and $y = x_l$.
A maximal $\Gamma$-connected subset of a given $B \subseteq G$
is said to be a \emph{$\Gamma$-connected component} of $B$.
\end{defn}

\begin{lem} \label{Folner_connected}
If $\varepsilon > 0$ and $\mu$ is a weighted $\varepsilon$-F\o lner set, then the support
of $\mu$ contains a $\Gamma$-connected component $A$
such that $\mu \restriction A$ is $\varepsilon$-F\o lner.
\end{lem}

\begin{proof}
Let $A_i$ $(i < l)$
enumerate the $\Gamma$-connected components of the support of $\mu$.
Then since
\[
 \sum_{i < l} \sum_{\gamma \in \Gamma} \sum_{s \in A_i} |\mu(s \cdot \gamma) - \mu(s)| <
\varepsilon  \sum_{i < l} \mu (A_i),
\]
there must exist an $i < l$ such that
\[
\sum_{\gamma \in \Gamma} \sum_{s \in A_i} |\mu(s \cdot \gamma) - \mu(s)| <
\varepsilon  \mu (A_i).
\]
Since $A_i$ is a $\Gamma$-connected component,
$\mu (s \cdot \gamma) = \mu \restriction A_i (s \cdot \gamma)$ whenever
$s$ is in $A_i$ and therefore $A_i$ is $\varepsilon$-F\o lner.
\end{proof}

This Lemma has the following useful consequence.
\begin{lem} \label{nice_Folner}
Let $G$ be a group with a finite generating set $\Gamma$, acting on itself
from the right.
If $A  \subseteq G$ is a $\varepsilon$-F\o lner set,
then there is a $B \subseteq G$ which is an $\varepsilon$-F\o lner set
such that $B$ is $\Gamma$-connected, $B$ contains the identity,
and $|B| \leq |A|$.
\end{lem}

\begin{proof}
Fix an $\varepsilon$-F\o lner set $A$ and
let $C$ be a $\Gamma$-connected component of $A$ which is
$\varepsilon$-F\o lner.
Let $g$ be any element of $C$ and define $B = g^{-1} C$.
It is easily verified that $C$ is still $\Gamma$-connected and
(right) $\varepsilon$-F\o lner.
\end{proof}

\section{F\o lner sets of trees}

\label{act_on_Tcal}

Rather than studying F\o lner sets in $F$ directly, it will
be easier to deal with weighted F\o lner sets in the partial right action of
$F$ on $\Tcal$ which I will
refer to as \emph{weighted F\o lner sets of trees}.
These are essentially weighted right F\o lner sets consisting of positive elements of $F$ (positive
with respect to the infinite presentation mentioned above).
The reformulation of the amenability problem for $F$ in terms of the existence of F\o lner
sets of positive elements is well known \cite{growth_amen_semi} and part of a more general phenomenon
(see \cite[1.28]{amenability}),
but we will need the more precise analytic consequences of the following lemmas.

\begin{lem} \label{not_in_marginal}
For every finite binary sequence $u$, the set $E_u$ of all $f \in F$ such that
$u$ is not extended by an element of $R_f$ is marginal with respect to the right action
of $F$ on itself.
\end{lem}

\begin{proof}
Let $u$ be fixed and
let $g$ be the element of $F$ which is the identity on
sequences which do not extend $u$ and which satisfies
\[
(u \tie v) \cdot g = u \tie (v \cdot x_0)
\]
for all sequences $v$ for which $v \cdot x_0$ is defined.

Now let $f$ be an element of $E_u$ and let $t$ be the initial part of $u$ which is
in $R_f$ (such a $t$ exists by our assumption that $u$ does not have an extension
in $R_f$).
Let $s = t \cdot f^{-1}$.
That is, $s$ is the element of $L_f$ such that
\[
|\{a \in L_f : a <_\lex s\}| = |\{b \in R_f : b <_\lex t\}|
\]
Observe that $f \cdot g^n$ can be represented by the tree diagram $(A,B)$ defined by:
\[
A = (L_f \setminus \{s\}) \cup  \{s \tie v : t \tie v \in L_{g^n}\}
\]
\[
B = (R_f \setminus \{t\}) \cup \{t \tie v : t \tie v \in R_{g^n}\}.
\]
That $(A,B)$ is a reduced tree diagram follows from the characterization mentioned in
Section \ref{background} and the following facts:
\begin{itemize}

\item $(L_f,R_f)$ and $(L_{g^n},R_{g^n})$ are reduced;

\item $|\{a \in A : a <_\lex s\}| = |\{b \in B : b <_\lex t\}|$;

\item
the minimum (maximum) elements of the sets
\[
\{s \tie v : t \tie v \in L_{g^n}\} \textrm{ and }
\{t \tie v : t \tie v \in R_{g^n}\}
\]
end in $\seq{0}$ (respectively $\seq{1}$).

\end{itemize}
Consequently $u$ is extended by an element of $R_{f \cdot g^n}$ for all $n > 0$ and
hence $g$ marginalizes $E_u$.
\end{proof}

\begin{lem} \label{Folner->trees}
There is a constant $C$ such that if $A \subseteq F$ is a (right)
$\varepsilon$-F\o lner set, then there is a weighted $C \varepsilon$-F\o lner set
of trees supported on a subset of $\{R_f : f \in A\}$.
\end{lem}

\begin{proof}
Let $U$ consist of all binary sequences of length 4, noting that
if $T$ dominates $U$, then every element of $\Gamma$ acts properly on $T$.
By Lemma \ref{not_in_marginal},
the set of $f \in F$ such that $U$ is not dominated by $R_f$
is marginal.
By Lemma \ref{marginal->Folner},
there is a constant $C > 1$ such that if $A$ is $\varepsilon$-F\o lner and
\[
A_0 = \{f \in A : U \textrm{ is dominated by } R_f\},
\]
then $A_0$ is $C \varepsilon$-F\o lner.
Observe that if $\gamma$ is a generator and $f$ is in $A_0$,
then $L_{f \cdot \gamma} = L_f$ and $R_{(f \cdot \gamma)} = (R_f) \cdot \gamma$.
We are now finished by Lemma \ref{Folner_hom_prop} applied to
$h(f) = R_f$ and $\mu  = 1_{A_0}$.
\end{proof}

\section{An operation on elements of $\Tcal$}

\label{partial_section}

In this section I will define an operation $\partial$
on elements of $\Tcal$ which reduces
their size logarithmically.

If $T$ is in $\Tcal$, then the \emph{end points} of $T$ are the maximum and
minimum elements of $T$.
All other elements of $T$ are said to be \emph{interior}.

\begin{defn}
Suppose that $T$ is in $\Tcal$.
$\partial T$ is the maximum $U \in \Tcal$ (with respect to the order
of domination) which is dominated by $T$
which satisfies the following \emph{defining conditions}:
\begin{enumerate}

\item \label{has_int} $U$ contains extensions of both $\seq{01}$ and $\seq{10}$;

\item \label{monotone_cond}
one of the following holds:
\begin{itemize}

\item 
if $u <_\lex v$ are interior elements of $U$,
then $2 |T/u| \leq |T/v|$;

\item 
if $u <_\lex v$ are interior elements of $U$,
then $2 |T/v| \leq |T/u|$;

\end{itemize}

\item \label{int_cond}
the minimum (respectively the maximum) interior element of $U$ terminates
with
a $\seq{1}$ (respectively with a $\seq{0}$).

\end{enumerate}
If no such $U$ exists, then $\partial T$ is defined to be the trivial tree.
\end{defn}

The following lemma justifies the use of \emph{maximum} in the definition of
$\partial T$.

\begin{lem}
If there is a $U$ satisfying the defining conditions for $\partial T$,
then there is a maximum such $U$ with respect to the order of domination.
\end{lem}

\begin{remark}
\emph{Condition \ref{int_cond}} is necessary to ensure the uniqueness
of maximal elements
of $\Tcal$ which satisfying the defining conditions for $\partial T$.
\end{remark}

\begin{proof}
Suppose for contradiction that $U$ and $V$ are distinct
maximal elements of $\Tcal$, each dominated by $T$,
which satisfy the defining conditions of $\partial T$.
First I claim that
\emph{condition \ref{monotone_cond}} is satisfied in the same way for $U$ and $V$.
For this, it is sufficient to show that if $W$ is in $\Tcal$, contains extensions
of both $\seq{01}$ and $\seq{10}$,
and satisfies the first (second) option of \emph{condition \ref{monotone_cond}}
for $T$, then 
$|T/\seq{01}| < |T/\seq{10}|$ (respectively $|T/\seq{10}| < |T/\seq{01}|$). 
Suppose that $W$ satisfies the first option (the other case is symmetric).
Let $w$ be the greatest element of $W$ extending $\seq{01}$ and $w'$ be the least element
of $W$ extending $\seq{10}$.
\emph{Condition \ref{monotone_cond}} implies that if $u_i$ $(i \leq l)$ are
the elements of $W$ which extend $\seq{01}$, then
\[
\sum_{i < l} |T/u_i| < |T/u_l| = |T/w| 
\]
and hence
$|T/\seq{01}| < 2|T/w|$.
It follows that
\[
|T/\seq{01}| < 2|T/w| \leq |T/w'| \leq |T/\seq{10}|.
\]

By replacing $U$ and $V$ with their mirror images if necessary,
we will assume that the quantity $|T/u|$ is increasing as
$u$ increases in $U$ (or equivalently as $u$ increases in $V$).
Since $U \symdif V$ is non-empty, there is a minimum element
of $U$ which either properly extends an element of $V$ or is properly extended
by an element of $V$.
By exchanging the roles of $U$ with $V$ if necessary, we may assume the former occurs.

First suppose that the minimum elements of $U$ and $V$ are the same.
Let $u$ be the greatest element of $U$ such that it and all of its
$<_\lex$-predecessors extend an element of $V$ (this
extension may not be proper in the case of the predecessors of $u$).
Let $v$ be the element of $V$ which $u$ extends.
Define
\[
W = \{w \in U : w \leq_\lex u\} \cup \{w \in V : v <_\lex w\}
\]
and observe that $W$ is in $\Tcal$.
Notice that $|T/u| \leq |T/v|$.
Also observe that
since $U$ does not dominate $V$, $u$ is not the maximum element of $U$ and
therefore if $x$ is an element of $U$ such that $x \leq_\lex u$, then
$|T/x| \leq |T/u|$.
It follows that if $x$ is an interior element of $W$ in $\{w \in U : w \leq_\lex u\}$ and $y$ is an interior element
of $W$ in $\{w \in V : v <_\lex w\}$,
then $x \leq_\lex u$ and $v <_\lex y$ which in turn implies
\[
2|T/x| \leq 2|T/u| \leq 2|T/v| \leq |T/y|.
\]
Since both $U$ and $V$ satisfy \emph{condition \ref{monotone_cond}}, 
if $x <_\lex y$ are interior elements of $W$ and either both are in $U$ or both are in $V$, then
$2|T/x| \leq |T/y|$.
It follows that
$W$  satisfies \emph{condition \ref{monotone_cond}} as well.
Since $U$ and $V$ have the same
minimum element and both satisfy \emph{condition \ref{int_cond}},
$U$ and $V$ have the same minimal interior element.
Thus the minimum interior element of $U$ is also the minimum interior element of $W$.
Also, since $V$ is not dominated by $U$, $u$ is not the maximum element of $W$ and $v$ is 
not the maximum element of $V$.
Moreover, it can not be the case that $v$ is the maximum interior element of $V$.
If this were the case, then $v = x \tie \seq{0}$ for some $x$.
It would then follow that
$x \tie \seq{1}$ would be the maximum element of $V$
and would have an extension in $U$.
This would imply that $U$ dominates $V$, which we assumed was not the case.
It follows that the maximum interior element of $W$ is the same as the maximum
interior element of $V$.
Therefore $W$ satisfies \emph{condition \ref{int_cond}}.
But now $W$ satisfies the defining conditions for $\partial T$,
contradicting the maximality of $V$.

Now suppose that the minimum elements of $U$ and $V$ differ.
Let $v$ be such that $v \tie \seq{0}$ is the minimum element of $V$, noting
that $v \tie \seq{1}$ is the minimum interior element of $V$.
Define
\[
W = \{w \in U : v \tie \seq{0} \subseteq w\} \cup \{w \in V : v \tie \seq{0} <_\lex w\}
\]
and observe that $W$ is in $\Tcal$.
Since $v \tie \seq{0}$ is not in $U$, it must be that both $v \tie \seq{00}$ and
$v \tie \seq{01}$ have extensions in $U$ and in particular, 
the minimum interior element of $W$ is the same as the one of $U$.
Since $V$ is not dominated by $U$, $v \ne \Seq{}$.
Since $v \tie \seq{0}$ is the minimum element of $V$, the entries of $v$ are
all $\seq{0}$ and hence
$v \tie \seq{1}$ is not the maximum element of
$V$ or of $W$.
Since every element of $W$ not extending $v \tie \seq{0}$ is in $V$, it follows
that the maximum interior element of $W$
is the same as the maximum interior element of $V$.
Therefore $W$ satisfies \emph{condition \ref{int_cond}}.
Next suppose that $x <_\lex y$ are interior elements of $W$.
If both $x$ and $y$ are in $U$ or both are in $V$, then
$2|T/x| \leq |T/y|$ follows from the fact that $U$ and $V$ satisfy \emph{condition \ref{monotone_cond}}
and that the minimum interior element of $W$ is the same as that of $U$ and that the maximum
interior element of $W$ is the same as that of $V$.
Next suppose that $x$ is in $\{w \in U : v \tie \seq{0} \subseteq w\}$ and $y$ is in $\{w \in V : v \tie \seq{0} <_\lex w\}$.
Observe that either $y = v \tie \seq{1}$ or $v \tie \seq{1} <_\lex y$.
Since $v \tie \seq{0}$ is extended by an element of $U$, it must be that $v \tie \seq{1}$ is also
extended by an element $w$ of $U$ which is in the interior of $U$.
We now have
\[
2|T/x| \leq |T/w| \leq |T/v\tie \seq{1}| \leq |T/y|
\]
(This is where the crucial use of \emph{condition \ref{int_cond}} occurs.)
Thus $W$ satisfies \emph{condition \ref{monotone_cond}}.
Again, $W$ satisfies the defining conditions for $\partial T$,
contradicting the maximality of $V$.
\end{proof}

\begin{lem} \label{partial_growth}
If $\partial T$ has $n$ elements, then $T$ has more than
$2^{n-2}$ elements.
\end{lem}

\begin{proof}
There are $n-2$ interior elements of $\partial T$ and
thus by \emph{condition \ref{monotone_cond}} the total
number of elements of $T$ which extend an interior element of $\partial T$
is at least $\sum_{i=1}^{n-2} 2^{i-1} = 2^{n-2} - 1$.
Since $\partial T$ is dominated by $T$, there are at least two elements of $T$ remaining
to be counted, putting the total greater than $2^{n-2}$.
\end{proof}

\begin{lem} \label{partial_compatible}
If $g$ is in $F$, $T$ is in $\Tcal$,
and $g$ acts properly on $\partial T$,
then $\partial (T \cdot g) = (\partial T) \cdot g$.
\end{lem}

\begin{proof}
First I will verify that if $g$ and $T$ are as in the statement of the lemma,
then $(\partial T) \cdot g$ is dominated by $\partial (T \cdot g)$.
Since the action of $g$ on the elements of $\partial T$
does not change their final digit,
$(\partial T) \cdot g$ satisfies \emph{condition \ref{int_cond}}.
Since the action of $g$ preserves lexicographic order and extension of sequences,
$(\partial T) \cdot g$ satisfies \emph{condition \ref{monotone_cond}}.
Finally, since the action of $g$ on $\partial T$ is proper,
$\partial T \cdot g$ satisfies \emph{condition \ref{has_int}}.
It follows that $(\partial T) \cdot g$ is dominated by $\partial (T \cdot g)$.

Next observe that if $g$ acts properly on $T$, then $g^{-1}$ acts
properly on $T \cdot g$.
It follows that
$\partial (T \cdot g) \cdot g^{-1}$ is dominated by 
$\partial (T \cdot g \cdot g^{-1}) = \partial T$.
Since $g$ and $g^{-1}$ are injections, it follows that $\partial (T \cdot g)$ and
$(\partial T) \cdot g$ have the same elements and hence are equal.
\end{proof}

\begin{defn}
Let $\Ecal$ be the set of all $T \in \Tcal$ such that neither of
the following inequalities hold:
\[
(+) \qquad |T/\seq{001}| < |T/\seq{01}| < |T/\seq{10}|
\]
\[
(-) \qquad |T/\seq{001}| > |T/\seq{01}| > |T/\seq{10}|
\]
Define
\[
\Tcal^+ = \{T \in \Tcal : |T/\seq{001}| < |T/\seq{01}| < |T/\seq{10}|\}
\]
\[
\Tcal^- = \{T \in \Tcal : |T/\seq{001}| > |T/\seq{01}| > |T/\seq{10}|\}
\]
\end{defn}

\begin{lem} \label{E_marginal}
$\Ecal$ is marginal.
\end{lem}

\begin{proof}

Define the following elements of $F$:
\[
a =
\begin{cases}
\seq{000} \mapsto \seq{000} \\
\seq{0010} \mapsto \seq{001} \\
\seq{0011} \mapsto \seq{0100} \\
\seq{01} \mapsto \seq{0101} \\
\seq{100} \mapsto \seq{011} \\
\seq{101} \mapsto \seq{10} \\
\seq{11} \mapsto \seq{11} \\
\end{cases}
b =
\begin{cases}
\seq{000} \mapsto \seq{000} \\
\seq{0010} \mapsto \seq{001} \\
\seq{0011} \mapsto \seq{01} \\
\seq{01} \mapsto \seq{100} \\
\seq{10} \mapsto \seq{101} \\
\seq{11} \mapsto \seq{11} \\
\end{cases}
\]
(i.e. $a = x_0^2 x_1 x_4 x_2^{-2} x_0^{-2}$ and $b = x_0^2 x_1 x_3^{-1} x_0^{-2}$).
Define
\[
\Ecal_a = \{T \in \Tcal : \max (|T/\seq{001}|,|T/\seq{10}|) = |T/\seq{01}| \}
\]
\[
\Ecal_b = \{T \in \Tcal : \max (|T/\seq{001}|,|T/\seq{10}|) < |T/\seq{01}| \}
\]
Notice that $T$ is in $\Ecal_{\max} = \Ecal_a \cup \Ecal_b$ if and only if
\[
|T/\seq{01}| = \max (|T/\seq{001}|,|T/\seq{01}|,|T/\seq{10}|).
\]
Furthermore, if $T$ is in $\Ecal \setminus \Ecal_{\max}$,
then
\[
|T/\seq{01}| = \min (|T/\seq{001}|,|T/\seq{01}|,|T/\seq{10}|).
\]
Observe that for all $T \in \Tcal$ such that $T \cdot a$ is defined, we have:
\[
|(T\cdot a)/\seq{001}| = |T/\seq{0010}| < |T/\seq{001}| 
\]
\[
|(T\cdot a)/\seq{01}| = |T/\seq{0011}| + |T/\seq{01}| + |T/\seq{100}| > |T/\seq{01}| 
\]
\[
|(T\cdot a)/\seq{10}| = |T/\seq{101}| < |T/\seq{10}| 
\]
Therefore if $T \cdot a$ is in $(\Ecal_a \cup \Ecal_b) \cdot a$, then
\[
\max( |(T \cdot a)/{\seq{001}}|, |(T\cdot a)/{\seq{10}}| ) < |(T\cdot a)/{\seq{01}}|
\]
and hence $T \cdot a$ is in $\Ecal_b$.
Since $\Ecal_a$ is disjoint from $\Ecal_b$,
this shows that $a$ marginalizes $\Ecal_a$.

Observe that for all $T \in \Tcal$ such that $T \cdot b$ is defined, we have:
\[
|(T\cdot b)/\seq{001}| = |T/\seq{0010}| < |T/\seq{001}| 
\]
\[
|(T\cdot b)/\seq{01}| = |T/\seq{0011}| < |T/\seq{001}| 
\]
\[
|(T\cdot b)/\seq{10}| = |T/\seq{01}| + |T/\seq{10}| > |T/\seq{10}| 
\]
Define
\[
\Rcal = \{T \in \Tcal : \max( |T/{\seq{001}}|, |T/{\seq{01}}| ) < |T/{\seq{10}}| \}.
\]
The above inequalities show that $(\Ecal_b \cup \Rcal) \cdot b \subseteq \Rcal$.
Since $\Ecal_b$ is disjoint from $\Rcal$,
this shows that $b$ marginalizes $\Ecal_b$.

Clearly the elements of $\Ecal \setminus \Ecal_{\max}$ lie in one
of the following sets:
\[
\Ecal_1 = \{T \in \Tcal : |T/{\seq{001}}| > |T/{\seq{01}}| = |T/{\seq{10}}|\}
\]
\[
\Ecal_2 = \{T \in \Tcal : |T/{\seq{01}}| < \min(|T/{\seq{001}}| , |T/{\seq{10}}|)\}
\]
\[
\Ecal_3 = \{T \in \Tcal : |T/{\seq{001}}| = |T/{\seq{01}}| < |T/{\seq{10}}|\}
\]
The proof will therefore be complete once it has been show that $x_0$
marginalizes $\Ecal_1 \cup \Ecal_2$ off $\Ecal_{\max}$ and that
$x_0$ marginalizes $\Ecal_3$ off $\Ecal_{\max}$.

If $T$ is in $\Tcal$, then
\[
|(T\cdot x_0)/\seq{01}| = |T/\seq{001}|
\]
\[
|(T \cdot x_0)/\seq{10}| = |T/\seq{01}|
\]
If $T$ is in $\Tcal^- \cup \Ecal_1 \cup \Ecal_2$, then $|T/\seq{001}| > |T/\seq{01}|$
and thus $|(T\cdot x_0)/\seq{01}| > |(T\cdot x_0)/\seq{10}|$.
This yields:
\[
(\Tcal^- \cup \Ecal_1 \cup \Ecal_2) \cdot x_0 \subseteq \Tcal^- \cup \Ecal_{\max}
\]
If $T$ is in $\Ecal_3$, then $|T/\seq{001}| = |T/\seq{01}|$ and
thus $|(T\cdot x_0)/\seq{01}| = |(T\cdot x_0)/\seq{10}|$.
This yields:
\[
\Ecal_3 \cdot x_0 \subseteq \Ecal_1 \cup \Ecal_{\max}
\]
This shows that $x_0$ marginalizes $\Ecal_1 \cup \Ecal_3$ off $\Ecal_{\max}$ and
that $x_0$ marginalizes $\Ecal_3$ off $\Ecal_{\max}$.
\end{proof}

\begin{defn}
Let $\Ecal^*$ be the set of all $T \in \Tcal$ such that neither of
the following inequalities hold:
\[
(2\times) \qquad 2 |T/\seq{001}| \leq |T/\seq{01}| \leq \frac{1}{2} |T/\seq{10}|
\]
\[
(\frac{1}{2} \times) \qquad \frac{1}{2} |T/\seq{001}| \geq |T/\seq{01}| \geq 2 |T/\seq{10}|
\]
\end{defn}

Observe that if $T$ is not in $\Ecal^*$ and $T$ contains extensions
of both $\seq{01}$ and $\seq{10}$, then $\partial T$ is not the trivial
tree since $\{\seq{00},\seq{01},\seq{10},\seq{11}\}$ is then an element of $\Tcal$
which satisfies the defining conditions for $\partial T$.

\begin{lem} \label{doubling_inv}
If $T$ satisfies $(+)$ and $\gamma$ is a generator,
then either $T \cdot \gamma$ is undefined,
$T \cdot \gamma$ is in $\Ecal$, or else $T \cdot \gamma$
satisfies $(+)$ (and similarly for $(-)$).
In particular, if $\Acal \subseteq \Tcal \setminus \Ecal^*$ is $\Gamma$-connected,
then either all elements of $\Acal$ satisfy $(2 \times)$ or all elements
of $\Acal$ satisfy $(\frac{1}{2} \times)$.
\end{lem}

\begin{proof}
This follows from the following equalities which hold whenever
the relevant action is defined:
\[
T/\seq{001} = (T \cdot x_0)/\seq{01}
\]
\[
T/\seq{01} = (T \cdot x_0)/\seq{10}
\]
\[
T/\seq{001} = (T \cdot x_1^{\pm})/\seq{001}
\]
\[
T/\seq{01} = (T \cdot x_1^{\pm})/\seq{01}.
\] 
\end{proof}

\begin{lem} \label{doubling}
$\Ecal^*$ is marginal.
\end{lem}

\begin{remark}
Lemmas \ref{doubling_inv} and \ref{doubling} have non trivial qualitative
consequences for $F$-invariant probability measures on $\Tcal$.
If $\mu$ is an $F$-invariant probability measure on $\Tcal$, then $\mu(\Ecal^*) = 0$.
Furthermore, $\mu(\Tcal^+ \symdif (\Tcal^+ \cdot \Gamma)) = 0$ and therefore if $\mu$ is
additionally ergodic, then it must be that $\mu$ assigns measure 1 either to the set of
elements of $\Tcal$ which satisfy $(2 \times)$ or else to the set of those which
satisfy $(\frac{1}{2} \times)$.
\end{remark}

\begin{proof}
Define the following elements of $F$:
\[
c =
\begin{cases}
\seq{00} \mapsto \seq{0} \\
\seq{01} \mapsto \seq{100} \\
\seq{10} \mapsto \seq{101} \\
\seq{11} \mapsto \seq{11} \\
\end{cases}
d =
\begin{cases}
\seq{000} \mapsto \seq{00} \\
\seq{001} \mapsto \seq{010} \\
\seq{01} \mapsto \seq{011} \\
\seq{1} \mapsto \seq{1} \\
\end{cases}
\]
(i.e. $c = x_0 x_1^{-1}$ and $d = x_0^2 x_1^{-1} x_0^{-1}$).
Define
\[
\Xcal = \{T \in \Tcal^+ : 2 |T/\seq{01}| \leq |T/\seq{10}|\}
\]
\[
\Ecal_4 = \{T \in \Tcal^+ : 2 |T/\seq{01}| > |T/\seq{10}| \}
\]
\[
\Ecal_5 = \{T \in \Tcal^+ : 2|T/\seq{001}| > |T/\seq{01}| \}.
\]
I first claim that $(\Xcal \cup \Ecal_4) \cdot c \subseteq \Xcal \cup \Ecal$.
To see this, suppose that $T$ is in $\Tcal^+$.
Then
\[
2|(T\cdot c)/\seq{01}| = 2 |T/\seq{001}| <
|T/\seq{01}|+|T/\seq{10}| = |(T \cdot c)/\seq{10}|.
\]
and hence if $T\cdot c$ is in $\Tcal^+$, it is in $\Xcal$.
Since
$T \cdot c$ is not in $\Tcal^-$, it is either in $\Tcal^+$ or in $\Ecal$.
This proves the claim.
Since $\Ecal_4$ is disjoint from $\Xcal$, 
it follows that $c$ marginalizes $\Ecal_4$ off $\Ecal$ and hence that
$\Ecal_4$ is marginal according to Lemma \ref{E_marginal}.
Also, by Lemma \ref{doubling_inv},
$\Ecal_5 \cdot x_0 \subseteq \Ecal_4 \cup \Ecal$ and therefore
$\Ecal_5$ is marginal because $x_0$ marginalizes $\Ecal_5$ off
$\Ecal_4 \cup \Ecal$.

Next define
\[
\Ycal = \{T \in \Tcal^- :  |T/\seq{01}| \geq 2|T/\seq{10}|\}
\]
\[
\Ecal_6 = \{T \in \Tcal^- : |T/\seq{01}| < 2|T/\seq{10}| \}
\]
\[
\Ecal_7 = \{T \in \Tcal^- : |T/\seq{001}| < 2|T/\seq{01}| \}.
\]
Arguing as above, $(\Ycal \cup \Ecal_6) \cdot d \subseteq \Ycal \cup \Ecal$ and hence
$d$ marginalizes $\Ecal_6$ off $\Ecal$.
Also $\Ecal_7 \cdot x_0 \subseteq \Ecal_6 \cup \Ecal$ and
consequently both $\Ecal_6$ and $\Ecal_7$ are marginal.
Since $\Ecal^* = \Ecal \cup \bigcup_{i=4}^7 \Ecal_i$, $\Ecal^*$ is marginal and
the proof is complete.
\end{proof}

\begin{lem} \label{proper_act}
The set
\[
\{T \in \Tcal : \Gamma \textrm{ does not act properly on }\partial T \}
\]
is marginal.
\end{lem}

\begin{proof}
First observe that
\[
\seq{1}^{i+1} \seq{0} \cdot x_0^{-i} = \seq{10}
\]
\[
\seq{1}^i\seq{0} \cdot x_0^{-i} = \seq{01}.
\]
Hence if $T$ is in $\Tcal$ and $T \cdot x_0^{-i}$ is defined,
then
\[
|T/\seq{1}^i \seq{0}| = |(T\cdot x_0^{-i})/\seq{01}|
\]
\[
|T/\seq{1}^{i+1} \seq{0}| = |(T\cdot x_0^{-i})/\seq{10}|.
\]
By Lemma \ref{doubling}, $\Ecal^*$ is marginal.
Also, it follows immediately from the definitions that
for each $f$ in $F$, the set of all $T \in \Tcal$ for which $T \cdot f$ is undefined
is marginalized off the emptyset by $f$.
By Remark \ref{marginal_rem}, it follows that
\[
\Ecal^{**} = \bigcup_{i = 0}^{16} \Ecal^* \cdot x_0^i \cup
\{T \in \Tcal : \exists i \leq 16\ (T \cdot x_0^{-i} \textrm{ is undefined})\}
\]
is marginal as well.
Observe that if $T$ is not in $\Ecal^{**}$, then
$\{T \cdot x_0^{-i} : 0 \leq i \leq 16\}$ is a $\Gamma$-connected subset of $\Tcal \setminus \Ecal^*$ and therefore
by Lemma \ref{doubling_inv}, one of the following two assertions holds:
\[
(\forall i < 16) \quad 2|T/\seq{1}^i\seq{0}| < |T/\seq{1}^{i+1}\seq{0}| 
\]
\[
(\forall i < 16) \quad 2|T/\seq{1}^{i+1}\seq{0}|  < |T/\seq{1}^i\seq{0}|.
\]
That is, the set of all $T$ which satisfy neither of these assertions is marginal.

Now let $U$ consist of all binary sequences of length 4,
noting that $U$ satisfies \emph{conditions \ref{has_int} and \ref{int_cond}}.
Furthermore if $T$ is an element of $\Tcal$ such that
$\partial T$ dominates $U$, then every element of $\Gamma$ acts properly on $\partial T$.
Set $R = \{\seq{1}^i\seq{0} : i < 14\} \cup \{\seq{1}^{15}\}$ and let
$g$ denote the element of $F$ defined by the tree diagram $(U,R)$.
If $U$ fails to satisfy \emph{condition \ref{monotone_cond}} with respect to $T$, then
$U \cdot g$ fails to satisfy \emph{condition \ref{monotone_cond}} with respect to $T \cdot g$ and, 
in particular, $T \cdot g$ must be in $\Ecal^{**}$.
Therefore the set of $T$ such that $U$ does not satisfy the defining conditions for
$\partial T$ is marginalized by $g$ off $\Ecal^{**}$.
\end{proof}

\begin{lem} \label{partial_marginalize}
There is a constant $C$ such that if $\mu$ is a weighted $\varepsilon$-F\o lner
set of trees and $C \varepsilon \leq 1$,
then there is a weighted $C \varepsilon$-F\o lner set of trees
which is supported on a subset of
\[
\{\partial T : (\mu(T) > 0) \land (\partial T \textrm{ is non trivial}) \land
(\Gamma \textrm{ acts properly on } \partial T)\}.
\]
\end{lem}

\begin{proof}
By Lemmas \ref{marginal->Folner} and \ref{proper_act},
there is a $C$ such that if $\mu$ is a weighted
$\varepsilon$-F\o lner set of trees and
\[
\Acal = \{T \in \Tcal : (\mu(T) > 0) \land (\Gamma \textrm{ acts properly on } \partial T)\},
\]
then $\mu \restriction \Acal$ is $C \varepsilon$-F\o lner.
Now let such a $\mu$ be given and define $\Acal$ as above.
By Lemma \ref{partial_compatible},
\[
(\partial T) \cdot \gamma = \partial (T \cdot \gamma)
\]
whenever $\gamma$ is a generator and $T$ is in $\Acal$.
Applying Lemma \ref{Folner_hom_prop} to
$\mu \restriction \Acal$ and $h = \partial$ gives the desired conclusion.
\end{proof}

Now we are ready to complete the proof of Theorem \ref{tower_growth}.
I will first prove the following claim.
\begin{claim}\label{*}
There exists a constant $K > 1$ such that
if $A \subseteq F$ is a $K^{-n}$-F\o lner set,
then $A$ contains an element with a tree diagram whose trees each contain
at least $\exp_{n} (0)$ elements.
\end{claim}

\begin{proof}
By Lemmas \ref{Folner->trees} and \ref{partial_marginalize}
there is a constant $K > 1$ such that:
\begin{enumerate}

\item \label{Folner_to_trees}
if $A \subseteq F$ is an $\varepsilon$-F\o lner set and $K \varepsilon \leq 1$,
then there is a weighted $K \varepsilon$-F\o lner set of trees
$\mu$ with support contained in $\{R_f:f \in A\}$;

\item \label{derivative_Folner}
if $\mu$ is a weighted $\varepsilon$-F\o lner set of trees and $K \varepsilon \leq 1$,
then there is a weighted $K \varepsilon$-F\o lner set of trees
$\nu$ which is supported on a subset of 
\[
\{\partial T : (\mu(T) > 0) \land (\partial T \textrm{ is non trivial}) \land (\Gamma \textrm{ acts properly on } T)\}
\]

\end{enumerate}
Thus by (\ref{Folner_to_trees}), if $A \subseteq F$ is $K^{-(n+1)}$-F\o lner, then
there is an $A' \subseteq A$ such that $\{R_f : f \in A'\}$ is an $K^{-n}$-F\o lner set.
By applying (\ref{derivative_Folner}) $n$ times and observing that
weighted F\o lner sets have non-empty supports,
we have that there is an $f \in A$ such that $\partial^{n} R_f$ is non trivial.
Let $k_i = |\partial^{n-i} R_f|$ and observe that by Lemma \ref{partial_growth},
$k_0 \geq 4$ and $k_{i+1} > 2^{k_i - 2}$.
It follows by induction that $\exp_i (0) + 2 < k_i$ and in particular
that $R_f$ contains at least $\exp_{n}(0)$ elements.
\end{proof}

By Theorem 1 and Proposition 2 of \cite{metric_Thompson}, if $f$ is in $F$,
then the distance from $f$ to the identity is at least $(k-2)/3$, where $k$ is the common cardinality of the
trees in the reduced tree diagram for $f$.
In particular, if $k$ is at least $3$ --- the minimum cardinality of a tree in any diagram representing a
non-identity element --- then the distance is at least $k/16$.
It is easily verified that for all $n >0$,
$\frac{1}{16} \exp_{4n}(0) \geq \exp_{n}(0)$.
If $K$ is a constant which satisfies the conclusion of Claim \ref{*}, then
define $C = K^4$.

I now claim that if $A$ is $C^{-n}$ F\o lner, then $|A| \geq \exp_n(0)$.
To see this, let $A \subseteq F$ be $C^{-n}$-F\o lner.
By Lemma \ref{nice_Folner}, there is a finite $A' \subseteq F$ which is $\Gamma$-connected, $C^{-n}$-F\o lner, and
satisfies $|A'| \leq |A|$.
Since $A'$ is $K^{-4n}$-F\o lner, our choice of $K$ implies that $A'$
has an element $a$ whose reduced tree diagram
contains trees with at least $\exp_{4n}(0)$ leaves.
It follows that the distance from $a$ to the identity is at least
$\frac{1}{16} \exp_{4n}(0) \geq \exp_n(0)$.
Since $A'$ is $\Gamma$-connected, it must contain at least $\exp_n(0)$ elements and thus
$|A| \geq |A'| \geq \exp_n(0)$, establishing Theorem \ref{tower_growth}.

\def\Dbar{\leavevmode\lower.6ex\hbox to 0pt{\hskip-.23ex \accent"16\hss}D}

\end{document}